\overfullrule=0pt
\centerline {\bf A characterization of the existence of zeros for operators with Lipschitzian derivative and closed range}\par
\bigskip
\bigskip
\centerline {BIAGIO RICCERI}
\bigskip
\bigskip
ABSTRACT. Let $H$ be a real Hilbert space and $\Phi:H\to H$ be a $C^1$ operator with Lipschitzian derivative and closed range. We prove that $\Phi^{-1}(0)\neq \emptyset$ if and only if,
for each $\epsilon>0$, there exist a convex set $X\subset H$ and a convex function $\psi:X\to {\bf R}$ such that
$\sup_{x\in X}(\|x\|^2+\psi(x))-\inf_{x\in X}(\|x\|^2+\psi(x))<\epsilon$ and $0\in \overline{\hbox {\rm conv}}(\Phi(X))$.\par
\bigskip
\bigskip
\bigskip
\bigskip
{\bf 1. Introduction and statements of the results}\par
\bigskip
In the sequel, $H$ is a real Hilbert space, with $\hbox {\rm dim}(H)\geq 2$, $\Omega$ is an open convex subset of $H$ and $\Phi:\Omega\to H$ is a given operator.\par
\smallskip
We say that $\Phi$ has a Lipschitzian derivative if $\Phi$ is Fr\'echet differentiable and the derivative of $\Phi$, denoted by $\Phi'$, is Lipschitzian, i.e.
one has
$$\sup_{x,y\in \Omega, x\neq y}{{\|\Phi'(x)-\Phi'(y)\|_{{\cal L}(H)}}\over {\|x-y\|}}<+\infty,$$
where ${\cal L}(H)$ is the space of all continuous linear operators from $H$ into itself endowed with the norm
$$\|T\|_{{\cal L}(H)}=\sup_{\|x\|\leq 1}\|T(x)\|.$$
For a generic set $A\subseteq H$, we denote by $\overline {\hbox {\rm conv}}(A)$ the closed convex hull of $A$, i.e. the smallest closed convex
set containing $A$.\par
\smallskip
We are interested in the existence of zeros for $\Phi$. \par
\smallskip
To shorten the statements, we now introduce the following notations. Namely, for each convex set $X\subseteq H$, we set
$$\delta_X:=\inf_{\psi\in \Gamma_X}\left(\sup_{x\in X}(\|x\|^2+\psi(x))-\inf_{x\in X}(\|x\|^2+\psi(x))\right),$$
where $\Gamma_X$ denotes the family of all convex functions $\psi:X\to {\bf R}$.
\smallskip
We have\par
\medskip
PROPOSITION 1.1. - {\it Let $X\subseteq H$ be a convex set with more than one point. Then, $\delta_X>0$.}\par
\medskip
We will get Proposition 1.1 as a by-product of Theorem 1.1 below. In [2], C. Zalinescu provided a proof of Proposition 1.1 based on quite
hard arguments of convex analysis.\par
\smallskip
Moreover, for each subset $V$ of $\Omega$, we denote by ${\cal A}_V$ the family of all convex sets $X\subseteq V$ such that
$$0\in \overline {\hbox {\rm conv}}(\Phi(X)).$$
The aim of this very short note is to establish the following result:\par
\medskip
THEOREM 1.1. - {\it Assume that $\Phi$ is a $C^1$ operator with Lipschitzian derivative and let $V$ be a subset of $\Omega$ such that
$\Phi(V)$ is closed. Then, the following assertions are equivalent:\par
\noindent
$(1)$\hskip 5pt $0\in\Phi(V)$;\par
\noindent
$(2)$\hskip 5pt  $\inf_{X\in {\cal A}_{V}}\delta_X=0$.}
\medskip
More precisely, the key result of this note is Theorem 1.2 below. Theorem 1.1 then follows as a by-product of it.\par
\medskip
THEOREM 1.2. - {\it Assume that $\Phi$ is a $C^1$ operator such that $\Phi'$ is Lipschitzian, with Lipschitz constant $L$. Let $V$ be a subset of
$\Omega$ such that
$$\eta:=\inf_{x\in V}\|\Phi(x)\|>0.$$
Then, for each convex set $X\subset V$  such that $\delta_X<{{2\eta}\over {L}}$,
one has
$$0\not\in\overline{\hbox {\rm conv}}(\Phi(X)).$$}\par
\medskip
We will prove Theorem 1.2 and Theorem 1.1 as well in the next section. The main tool that we will use is the classical Kneser minimax theorem
[1]. For the reader convenience, we recall its statement:\par
\medskip
THEOREM 1.A. - {\it Let $E$ be a vector space, $F$ a locally convex topological vector space, $X\subseteq E$ a convex set and $Y\subset F$ a compact convex set. Moreover, let $h:X\to {\bf R}$ be a function which is convex in $X$ and upper semicontinuous and concave in $Y$.
Then, one has
$$\sup_Y\inf_Xh=\inf_X\sup_Yh.$$}\par
\medskip
The interest of Theorem 1.1 resides mainly in its full novelty which is, in turn, due to the particular proof approach based on Theorem 1.A. 
Actually, we are not aware of known results with which Theorem 1.1 can be compared even in a vague manner.\par
\smallskip
 In particular, it seems that the number $\delta_X$ is here introduced for the first time. We think that it is worth studying in depth. For instance, it would be useful to provide
an explicit positive lower bound for it. Indeed, the current proofs of Proposition 1.1 (ours and that in [2]) are highly indirect.\par
\smallskip
Another interesting feature of Theorem 1.1 resides in the derivative of $\Phi$: not only it does not appear in the conclusion at all, but also it
needs to be Lipschitzian. In this connection, the following example is enlightening.\par
\medskip
EXAMPLE 1.1. - Let $h:{\bf R}\to {\bf R}$ be any $C^1$ function with the following property: $[-\pi,\pi]\subseteq h({\bf R})$ and there exist two sequences in $(0,+\infty)$ $\{\alpha_n\}, \{\beta_n\}$ such that 
$$\lim_{n\to \infty}\alpha_n=+\infty,$$
  $$\lim_{n\to \infty}(\beta_n^2-\alpha_n^2)=0,$$
$$\alpha_n<\beta_n,$$ $$h(\alpha_n)=0,$$ $$h(\beta_n)=-\pi$$ for all $n\in {\bf N}$. 
For instance, one can take: $h(x)=\pi\sin(2\pi x^2)$, $\alpha_n=\sqrt{n}$, $\beta_n=\sqrt{n+{{3}\over {4}}}$.
Then, consider the function $\Phi:{\bf R}^2\to {\bf R}^2$ defined by
$$\Phi(x,y)=(\sin(h(x)), \cos(h(x)))$$
for all $(x,y)\in {\bf R}^2$. So, $\Phi$ is $C^1$ and
$$\Phi({\bf R}^2)=\{(s,t)\in {\bf R}^2 : s^2+t^2=1\}.$$
In connection with Theorem 1.1, take
$\Omega=V={\bf R}^2$. Now, fix $\epsilon>0$ and $n\in {\bf N}$ so that $\beta_n^2-\alpha_n^2<\epsilon$. Set
$$X:=[\alpha_n,\beta_n]\times \{0\}.$$
So, $X$ is convex and the points $(0,1)$, $(0,-1)$ belong to $\Phi(X)$. Consequently
$$0\in \hbox {\rm conv}(\Phi(X)).$$
Therefore, $X\in {\cal A}_V$. Moreover, we have
$$\delta_X\leq \sup_{(x,y)\in X}\|(x,y)\|^2-\inf_{(x,y)\in X}\|(x,y)\|^2=\beta_n^2-\alpha_n^2<\epsilon.$$
This shows that
$$\inf_{Y\in {\cal A}_V}\delta_Y=0.$$
However, $\Phi(V)$ is closed and $0\not\in\Phi(V)$. So, the implication $(2)\to (1)$ in Theorem 1.1 does not hold.
\smallskip

\bigskip
{\bf 2. Proofs}\par
\bigskip
{\it Proof of Theorem 1.2}. Fix any convex set $X\subset V$ and any convex function $\psi:X\to {\bf R}$ satisfying
$$\sup_{x\in X}(\|x\|^2+\psi(x))-\inf_{x\in X}(\|x\|^2+\psi(x))<{{2\eta}\over {L}}.\eqno{(2.1)}$$
Of course, pairs $(X,\psi)$ of this type do exist: for instance, this is the case when $X$ is a singleton.
Set
$$Y=\{x\in H : \|x\|\leq 1\}$$
and consider the functions $\varphi:X\to {\bf R}$, $f:\Omega\times Y\to {\bf R}$ and $g:X\times Y\to {\bf R}$ defined by
$$\varphi(x)={{L}\over {2}}(\|x\|^2+\psi(x)),$$
$$f(x,y)=\langle \Phi(x),y\rangle,$$
and
$$g(x,y)=f(x,y)+\varphi(x).$$
We claim that
$$\inf_X\sup_Yf-\sup_Y\inf_Xf\leq \sup_X\varphi-\inf_X\varphi.\eqno{(2.2)}$$
Arguing by contradiction, assume that
$$\inf_X\sup_Yf-\sup_Y\inf_Xf> \sup_X\varphi-\inf_X\varphi.$$
We then would have
$$\sup_Y\inf_Xg\leq \sup_Y\inf_Xf+\sup_X\varphi<\inf_X\sup_Yf+\inf_X\varphi\leq
\inf_X\sup_Yg. \eqno{(2.3)}$$
 For each $y\in Y$, the function $f(\cdot,y)$ is $C^1$ 
and its derivative (denoted by $f'_x(\cdot,y)$) is given by
$$\langle f'_x(x,y),u\rangle=\langle\Phi'(x)(u),y\rangle$$
for all $x\in \Omega$, $u\in H$, Also, for each $v, w\in\Omega$, we have
$$\|f'_x(v,y)-f'_x(w,y)\|=\sup_{u\in Y}|\langle\Phi'(v)(u)-\Phi'(w)(u),y\rangle|\leq
\sup_{u\in Y}\|\Phi'(v)(u)-\Phi'(w)(u)\|$$
$$=\|\Phi'(v)-\Phi'(w)\|_{{\cal L}(H)}\leq
L\|v-w\|.\eqno{(2.4)}$$
Taking $(2.4)$ into account, we also have
$$\langle Lv+f'_x(v,y)-Lw-f'_x(w,y),v-w\rangle$$
$$\geq L\|v-w\|^2-\|f'_x(v,y)-f'_x(w,y)\|\|v-w\|\geq L\|v-w\|^2-L\|v-w\|^2=0.$$
In other words, the derivative of the function ${{L}\over {2}}\|\cdot\|^2+f(\cdot,y)$ is monotone in $\Omega$ and so the function is
convex there ([3], Proposition 42.6). This implies that
$g(\cdot,y)$ is convex in $X$ since $\psi$ is so. Furthermore, $Y$ is weakly compact and, for each $x\in X$, the function $g(x,\cdot)$ is weakly
continuous, being affine and continuous. Hence, applying Theorem 1.A, we would have
$$\sup_Y\inf_Xg=\inf_X\sup_Yg$$
contradicting $(2.3)$. So, $(2.2)$ does hold. Notice that
$$\inf_X\sup_Yf=\inf_{x\in X}\|\Phi(x)\|.$$
Therefore, from $(2.2)$ we infer that
$${{L}\over {2}}\left(\inf_{x\in X}(\|x\|^2+\psi(x))-\sup_{x\in X}(\|x\|^2+\psi(x))\right)+\inf_{x\in X}\|\Phi(x)\|\leq \sup_{y\in Y}\inf_{x\in X}\langle \Phi(x),y\rangle$$
and hence, in view of $(2.1)$, taking into account that $\eta\leq \inf_{x\in X}\|\Phi(x)\|$, we have
$$0<\sup_{y\in Y}\inf_{x\in X}\langle \Phi(x),y\rangle.$$
Because of this inequality, we can fix $\gamma>0$ and $\tilde y\in Y$ so that
$$\inf_{x\in X}\langle\Phi(x),\tilde y\rangle\geq \gamma.$$
Thus, if we set
$$C=\{u\in H : \langle u,\tilde y\rangle\geq \gamma\},$$
we have
$$\overline {\hbox {\rm conv}}(\Phi(X))\subset C$$
since $C$ is closed and convex,
while $0\not\in C$ and the proof is complete.\hfill $\bigtriangleup$
\medskip
{\it Proof of Theorem 1.1.} The implication $(1)\to (2)$ is immediate. Indeed, if there exists $\tilde x\in V$ such that $\Phi(\tilde x)=0$, then
the singleton $\{\tilde x\}$ belongs to the family ${\cal A}_{V}$ and $\delta_{\{\tilde x\}}=0$, and so $(2)$ holds.
Viceversa, assume that $(2)$ holds. We have to prove $(1)$. Arguing by
contradiction, suppose that $0\not\in \Phi(V)$. Since $\Phi(V)$ is closed, this implies that
$\inf_{x\in V}\|\Phi(x)\|>0$. By assumption, there is a convex set $X\subset V$ 
such that
$$\delta_X<2{{\inf_{x\in V}\|\Phi(x)\|}\over {L}}$$
and
$$0\in \overline {\hbox {\rm conv}}(\Phi(X)),$$
contradicting Theorem 1.2.\hfill $\bigtriangleup$
\medskip
From the proof of Theorem 1.1, it clearly follows the following\par
\medskip
THEOREM 2.1. - {\it Assume that $\Phi$ is a $C^1$ operator with Lipschitzian derivative and let $V$ be a subset of $\Omega$ such that
$$\inf_{X\in {\cal A}_{V}}\delta_X=0.$$
Then, $\inf_{x\in V}\|\Phi(x)\|=0$.}\par
\medskip
{\it Proof of Proposition 1.1.} Arguing by contradiction, assume that there is a convex set $X\subset H$, with more than one point, such that
$\delta_X=0$. Fix $x_1, x_2\in X$, with $x_1\neq x_2$. Let $V$ be the closed segment joining $x_1$ with $x_2$. Clearly, $\delta_V=0$.
Fix $u, v, w\in H$ so that
$\langle u,v\rangle =0$, $\|u\|=\|v\|=1$, $\langle w,x_1\rangle\neq \langle w,x_2\rangle$. Finally, consider the operator
$\Phi:H\to H$ defined by
$$\Phi(x)=\sin\left({{\langle w,x-x_1\rangle\pi}\over {\langle w,x_2-x_1\rangle}}\right)u+\cos\left({{\langle w,x-x_1\rangle\pi}\over {\langle w,x_2-x_1\rangle}}\right)v.$$
Of course, $\Phi$ has a Lipschitzian derivative, $\Phi(V)$ is compact and $0\not\in \Phi(V)$. Moreover, $\Phi(x_1)=v$ and $\Phi(x_2)=-v$.
Consequently, $0\in {\rm conv}(\Phi(V))$ and so $V\in {\cal A}_V$. Hence, condition $(2)$ of Theorem 1.1 is satisfied and not condition $(1)$. This
contradiction ends the proof.\hfill $\bigtriangleup$
\bigskip
\bigskip
\bigskip
\bigskip
{\bf Acknowledgements:} This work has been funded by the European Union - NextGenerationEU Mission 4 - Component 2 - Investment 1.1 under the Italian Ministry of University and Research (MUR) programme "PRIN 2022" - grant number 2022BCFHN2 - Advanced theoretical aspects in PDEs and their applications - CUP: E53D23005650006. The author has also been supported by the Gruppo Nazionale per l'Analisi Matematica, la Probabilit\`a e 
le loro Applicazioni (GNAMPA) of the Istituto Nazionale di Alta Matematica (INdAM) and by the Universit\`a degli Studi di Catania, PIACERI 2020-2022, Linea di intervento 2, Progetto ”MAFANE”. The author wishes to thank the referee for his/her remarks showing the excessive conciseness of the original version.\par
\bigskip
\bigskip
\bigskip
\bigskip
\centerline {\bf References}\par
\bigskip
\bigskip
\noindent
[1]\hskip 5pt  H. KNESER, {\it Sur un th\'eor\`eme fondamental de la th\'eorie des jeux}, C. R. Acad. Sci. Paris {\bf 234} (1952), 2418-2420.\par
\smallskip
\noindent
[2]\hskip 5pt C. ZALINESCU, personal communication.\par
\smallskip
\noindent
[3]\hskip 5pt E. ZEIDLER, {\it Nonlinear functional analysis and its
applications}, vol. III, Springer-Verlag, 1985.\par
\bigskip
\bigskip
\bigskip
\bigskip
Department of Mathematics and Informatics\par
University of Catania\par
Viale A. Doria 6\par
95125 Catania, Italy\par
{\it e-mail address}: ricceri@dmi.unict.it

\bye